\documentclass[12pt]{article}
\usepackage{epsfig,amssymb}
\usepackage[active]{srcltx}
\usepackage{rotating,booktabs}

\begin{document}
\newcommand{\Gm}{G_\mu}
\newcommand{\Pm}{P_\mu}
\newcommand{\Pim}{\Pi_\mu}
\newcommand{\Hm}{H_\mu}
\newcommand{\Lm}{L_\mu}
\newcommand{\Thm}{\Xi_\mu}
\newcommand{\Xim}{\Xi_\mu}

\newcommand{\ux}{u_{,x}}
\newcommand{\uxx}{u_{,xx}}
\newcommand{\ut}{u_{,t}}
\newcommand{\utt}{u_{,tt}}
\newcommand{\vt}{v_{,t}}
\newcommand{\ppt}{p_{,t}}
\newcommand{\ob}{b}
\newcommand{\cc}{c}
\newcommand{\sh}{\sharp}
\newcommand{\fl}{\flat}
\newcommand{\odeg}{\; ; \;}
\newcommand{\bif}{\; , \;}
\newcommand{\bis}{\; ; \;}
\newcommand{\bit}{\; , \;}

\newcommand{\hfillbox}{\hfill $\Box$}
\newcommand{\Gf}{{\Gamma_f}}
\newcommand{\sta}{{\rm sta}}
\newcommand{\ext}{{\rm ext}}
\newcommand{\extu}{{\ext}_{\hspace{-.5cm}{}_{{}_{{}_{u\in\calU}}}}}
\newcommand{\extx}{{\ext}_{\hspace{-.5cm}{}_{{}_{{}_{\x\in\calX}}}}}
\newcommand{\exty}{{\ext}_{\hspace{-.5cm}{}_{{}_{{}_{\y\in\calY}}}}}
\newcommand{\extyd}{{\ext}_{\hspace{-.5cm}{}_{{}_{{}_{\y^*\in\calY*}}}}}
\newcommand{\extud}{{\ext}_{\hspace{-.5cm}{}_{{}_{{}_{u^*\in\calU^*}}}}}
\newcommand{\extua}{{\ext}_{\hspace{-.5cm}{}_{{}_{{}_{u\in\calU_a}}}}}
\newcommand{\exte}{{\ext}_{\hspace{-.5cm}{}_{{}_{{}_{\eps\in\calE}}}}}
\newcommand{\extee}{{\ext}_{\hspace{-.5cm}{}_{{}_{{}_{e\in\calE}}}}}
\newcommand{\exted}{{\ext}_{\hspace{-.5cm}{}_{{}_{{}_{\eps^*\in\calE^*}}}}}
\newcommand{\extxn}{{\ext}_{\hspace{-.5cm}{}_{{}_{{}_{\bx\in\real^n}}}}}
\newcommand{\staua}{{\sta}_{\hspace{-.5cm}{}_{{}_{{}_{u\in\calU_a}}}}}
\newcommand{\stax}{{\sta}_{\hspace{-.5cm}{}_{{}_{{}_{\x\in\calX}}}}}
\newcommand{\staxa}{{\sta}_{\hspace{-.5cm}{}_{{}_{{}_{\x\in\calX_a}}}}}
\newcommand{\stayd}{{\sta}_{\hspace{-.5cm}{}_{{}_{{}_{\y^*\in\calY^*}}}}}
\newcommand{\stayad}{{\sta}_{\hspace{-.5cm}{}_{{}_{{}_{\y^*\in\calY^*_a}}}}}
\newcommand{\bEn}{{\bf E}^{(\eta)}}
\newcommand{\bTn}{{\bf T}^{(\eta)}}

\newcommand{\calTn}{\cal{T}^{(\eta)}}
\newcommand{\calEn}{\cal{E}^{(\eta)}}

\newcommand{\calA}{{\cal{A}}}
\newcommand{\barbA}{\bar{\bf A}}
\newcommand{\calO}{{\cal{O}}}
\newcommand{\calXo}{{\cal{X}}^o}
\newcommand{\calXk}{{\cal{X}}_k}

\newcommand{\U}{U}
\newcommand{\Ub}{\bar{U}}
\newcommand{\Uv}{\check{U}}
\newcommand{\Un}{\hat{U}}
\newcommand{\Ut}{\tilde{U}}
\newcommand{\Vv}{\check{V}}
\newcommand{\Vc}{V^c}
\newcommand{\Vn}{\hat{V}}
\newcommand{\Vb}{\bar{V}}
\newcommand{\Vt}{\tilde{V}}
\newcommand{\Gamv}{\check{\Gamma}}
\newcommand{\Gamn}{\hat{\Gamma}}
\newcommand{\Wb}{\bar{W}}
\newcommand{\Pb}{\bar{P}}
\newcommand{\Fb}{\bar{F}}

\newcommand{\cof}{{\rm cof}}
\newcommand{\adj}{{\rm adj}}
\newcommand{\meas}{{\rm meas}}
\newcommand{ \Ii }{{\partial I}}
\newcommand{ \Iu }{{\partial I_u}}
\newcommand{ \It }{{\partial I_f}}
\newcommand{\barIu}{{\bar{I}_u}}
\newcommand{\barIt}{{\bar{I}_f}}
\newcommand{ \Ji }{J_{\partial I}}
\newcommand{ \Ju }{J_{u}}
\newcommand{ \Js }{J_{\sig}}
\newcommand{ \Jj }{{\cal{J}}}
\newcommand{ \calI }{{\cal{I}}}
\newcommand{ \HH}{{H}}
\newcommand{ \PP}{{P}}
\newcommand{ \LL}{{L}}
\newcommand{ \JJ}{{J}}
\newcommand{ \II}{{\Psi}}
\newcommand{\calLs}{{L_{\mbox{\sta}}}}

\newcommand{ \id }{{i_d}}
\newcommand{\bell}{{\mbox{\large $\ell$}}}
\newcommand{\vsig}{\varsigma}
\newcommand{\vrho}{\varrho}
\newcommand{\vsigbar}{\bar{\varsigma}}
\newcommand{\barvsig}{\bar{\varsigma}}
\newcommand{\vsigb}{\bar{\varsigma}}

\newcommand{\pp}{p}
\newcommand{\vv}{v}
\newcommand{\Ww}{\bar{W}}
\newcommand{\ww}{w}
\newcommand{\barww}{\bar{w}}
\newcommand{\Ll}{{\Large \ell}}
\newcommand{\barbw}{\bar{\bf w}}
\newcommand{\uu}{u}
\newcommand{\zz}{z}
\newcommand{\barpp}{\bar{p}}
\newcommand{\barss}{\bar{s}}
\newcommand{\barx}{\bar{x}}
\newcommand{\bary}{\bar{y}}
\newcommand{\qq}{q}
\newcommand{\barqq}{\bar{\qq}}

\newcommand{\jj}{\mbox{\large{$\jmath$}}}
\newcommand{\cl}{\mbox{cl}}

\newcommand{\Lamo}{{\Lam^o}}
\newcommand{\Bo}{{B^o}}
\newcommand{\pol}{\; ; \; }
\newcommand{\oLam}{\Lamo^*}
\newcommand{\calUo}{{\calU^o}}
\newcommand{\Ao}{{A^o}}
\newcommand{\uo}{{u^o}}
\newcommand{\baruo}{{\bar{u}^o}}
\newcommand{\barou}{{\bar{u}^{o*}}}

\newcommand{\Co}{{C^o}}
\newcommand{\ocalX}{\calX^{o*}}
\newcommand{\ocalU}{\calUo^*}
\newcommand{\ou}{u^{o*}}
\newcommand{\infc}{\stackrel{+}{\vee}}
\newcommand{\rarw}{\rightarrow}
\newcommand{\rarwr}{\rightarrow \real}
\newcommand{\rawu}{\rightharpoonup}
\newcommand{\rawuast}{\stackrel{\ast}{\rightharpoonup}}
\newcommand{\pl}{{\parallel}}
\newcommand{\del}{{\delta}}
\newcommand{\barpartial}{{\bar{\partial}}}
\newcommand{\real}{{\mathbb R}} 
\newcommand{\bolM}{{\mathbb M}} 
\newcommand{\boI}{{\mathbb I}}
\newcommand{\boJ}{{\mathbb J}}

\newcommand{\comp}{{\mathbb C}}
\newcommand{\field}{{\mathbb F}}
\newcommand{\bareal}{\bar{\real}}
\newcommand{\realr}{\stackrel{\rightharpoonup}{\real}}
\newcommand{\reall}{\stackrel{\leftharpoonup}{\real}}
\newcommand{\rareal}{{\rightarrow \real}}
\newcommand{\rabareal}{{\rightarrow \bareal}}
\newcommand{\kap}{{\kappa}}
\newcommand{\dotp}{\dot{p}}
\newcommand{\dotkap}{\dot{\kappa}}
\newcommand{\dotkappa}{\dot{\kappa}}
\newcommand{ \doteps}{\dot{\epsilon}}
\newcommand{ \dottau}{\dot{\tau}}
\newcommand{ \doteta}{\dot{\eta}}
\newcommand{ \dotzet}{\dot{\zeta}}
\newcommand{\dotepsilon}{\dot{\epsilon}}
\newcommand{\half}{\frac{1}{2}}
\newcommand{\Ss}{{\bf \Sigma }}
\newcommand{\biota}{\mbox{\boldmath$\iota$}}
\newcommand{\syst}{{\mathbb S}} 
\newcommand{\systo}{{\syst^o}}
\newcommand{\bsig}{\mbox{\boldmath$\sigma$}}
\newcommand{\bvsig}{\mbox{\boldmath$\varsigma$}}
\newcommand{\bmu}{\mbox{\boldmath$\mu$}}
\newcommand{\barbmu}{\bar{\bmu}}
\newcommand{\brho}{\mbox{\boldmath$\rho$}}
\newcommand{\veps}{\xi}
\newcommand{\barveps}{\bar{\xi}}
\newcommand{\bveps}{\mbox{\boldmath$\xi$}}
\newcommand{\bfeta}{\mbox{\boldmath$\eta$}}
\newcommand{\barbsig}{\bar{\mbox{\boldmath$\sigma$}}}
\newcommand{\barbvsig}{\bar{\mbox{\boldmath$\varsigma$}}}
\newcommand{\barbpsi}{\bar{\mbox{\boldmath$\psi$}}}
\newcommand{\beps}{\mbox{\boldmath$\epsilon$}}
\newcommand{\barbeps}{\bar{\mbox{\boldmath$\epsilon$}}}

\newcommand{\bPhi}{{\mbox{\boldmath$\Phi$}}}
\newcommand{\bPsi}{{\mbox{\boldmath$\Psi$}}}
\newcommand{\btau}{{\mbox{\boldmath$\tau$}}}
\newcommand{\bbta}{{\mbox{\boldmath$\beta$}}}
\newcommand{\blam}{{\mbox{\boldmath$\lambda$}}}
\newcommand{\bLam}{{\mbox{\boldmath$\Lambda$}}}
\newcommand{\barbtau}{\bar{\mbox{\boldmath$\tau$}}}
\newcommand{\bchi}{{\mbox{\boldmath$\chi$}}}
\newcommand{\bchio}{\bchi^o}
\newcommand{\barbchio}{\bar{\bchi}^o}
\newcommand{\barbchi}{\bar{\mbox{\boldmath$\chi$}}}
\newcommand{\bxi}{{\mbox{\boldmath$\xi$}}}
\newcommand{\bphi}{{\mbox{\boldmath$\phi$}}}
\newcommand{\barphi}{\bar{\phi}}
\newcommand{\barbphi}{\bar{\mbox{\boldmath$\phi$}}}
\newcommand{\bpsi}{{\mbox{\boldmath$\psi$}}}
\newcommand{\barbxi}{\bar{\mbox{\boldmath$\xi$}}}
\newcommand{\bnu}{{\mbox{\boldmath$\nu$}}}
\newcommand{\barI}{{\bar{I}}}
\newcommand{\barK}{{\bar{K}}}
\newcommand{\barU}{{\bar{U}}}
\newcommand{\barV}{{\bar{V}}}
\newcommand{\hatV}{{\hat{V}}}
\newcommand{\barW}{{\bar{W}}}
\newcommand{\hatW}{{\hat{W}}}
\newcommand{\hatau}{{\hat{\tau}}}
\newcommand{\hatsig}{{\hat{\sigma}}}
\newcommand{\barv}{{\bar{v}}}
\newcommand{\baru}{\bar{u}}
\newcommand{\barh}{\bar{h}}
\newcommand{\bars}{\bar{s}}
\newcommand{\barp}{\bar{p}}
\newcommand{\barq}{\bar{q}}
\newcommand{\la}{\langle}
\newcommand{\ra}{\rangle}
\newcommand{\barho}{\bar{\rho}}
\newcommand{\bartheta}{\bar{\theta}}
\newcommand{\barbe}{\bar{\bf e}}
\newcommand{\barbeta}{\bar{\beta}}
\newcommand{\bareps}{\bar{\epsilon}}
\newcommand{\barmu}{\bar{\mu}}

\newcommand{\barsig}{\bar{\sigma}}
\newcommand{\barchi}{\bar{\chi}}
\newcommand{\barM}{\bar{M}}
\newcommand{\barxi}{\bar{\xi}}
\newcommand{\bartau}{\bar{\tau}}
\newcommand{\barpsi}{\bar{\psi}}
\newcommand{\ba}{{{\bf a}}}
\newcommand{\br}{{{\bf r}}}
\newcommand{\bp}{{{\bf p}}}
\newcommand{\Oo}{\Omega}
\newcommand{\Oot}{{\Omega_t}}
\newcommand{\Gxi}{{\Gamma_\chi}}
\newcommand{\Gt}{{\Gamma_t}}
\newcommand{\Gu}{{\Gamma_u}}
\newcommand{\bE}{{\bf E}}
\newcommand{\bQ}{{\bf Q}}
\newcommand{\ep}{{\epsilon}}
\newcommand{\bM}{{\bf M}}
\newcommand{\bK}{{\bf K}}
\newcommand{\eba}{\begin{array}}
\newcommand{\eea}{\end{array}}
\newcommand{\ebe}{\begin{eqnarray}}
\newcommand{\eee}{\end{eqnarray}}
\newcommand{\eb}{\begin{equation}}
\newcommand{\ee}{\end{equation}}
\newcommand{\bT}{{\bf T}}
\newcommand{\bJ}{{\bf J}}
\newcommand{\barcalU}{{\bar{\cal{U}}}}
\newcommand{\barcalE}{{\bar{\cal{E}}}}
\newcommand{\pcalU}{{\bar{\cal{U}}}}
\newcommand{\pcalE}{{\bar{\cal{E}}}}
\newcommand{\calW}{{\cal{W}}}
\newcommand{\calP}{{\cal{P}}}
\newcommand{\calPma}{\calP_{\max}}
\newcommand{\calPmi}{\calP_{\min}}
\newcommand{\calPst}{\calP_{\rm sta}}
\newcommand{\calPi}{\calP_{\inf}}
\newcommand{\calPs}{\calP_{\sup}}
\newcommand{\calPvi}{\calP_{\rm vi}}
\newcommand{\calPe}{\calP_{\rm ext}}
\newcommand{\calPbv}{\calP_{\rm bv}}
\newcommand{\calPcbv}{\calP_{\rm cbv}}
\newcommand{\calPiv}{\calP_{\rm iv}}
\newcommand{\calLe}{L_{\rm ext}}
\newcommand{\calPm}{\calP_{\min}}
\newcommand{\Pbv}{\calP_{\rm bv}}
\newcommand{\calPbvc}{\calP_{\rm cbv}}
\newcommand{\Bvp}{\calP_{\rm vp}}

\newcommand{\sT}{\cal{T}}
\newcommand{\dsT}{\dot{\cal{T}}}
\newcommand{\bg}{{\bf g}}
\newcommand{\bm}{{\bf m}}
\newcommand{\bC}{{\bf C}}
\newcommand{\bG}{{\bf G}}
\newcommand{\bH}{{\bf H}}
\newcommand{\bn}{{\bf n}}
\newcommand{\bt}{{\bf t}}
\newcommand{\bff}{{\bf f}}
\newcommand{\bV}{{\bf V}}
\newcommand{\bY}{{\bf Y}}
\newcommand{\BR}{{\bf R}}
\newcommand{\bR}{{\bf R}}
\newcommand{\bv}{{\bf v}}
\newcommand{\bw}{{\bf w}}
\newcommand{\BN}{{\bf N}}
\newcommand{\bx}{{\bf x}}
\newcommand{\by}{{\bf y}}
\newcommand{\BX}{{\bf X}}
\newcommand{\bX}{{\bf X}}
\newcommand{\BF}{{\bf F}}
\newcommand{\bN}{{\bf N}}
\newcommand{\bF}{{\bf F}}
\newcommand{\bA}{{\bf A}}
\newcommand{\bB}{{\bf B}}
\newcommand{\bI}{{\bf I}}
\newcommand{\bj}{{\bf j}}
\newcommand{\bk}{{\bf k}}
\newcommand{\calS}{{\cal S}}
\newcommand{\calSs}{{\calS_s}}
\newcommand{\calSa}{{\calS_a}}
\newcommand{\calM}{{\cal M}}
\newcommand{\calMs}{{\calM_{sym}}}
\newcommand{\calMp}{\calM_+}
\newcommand{\calMo}{\calM_{ort}}
\newcommand{\DelM}{{\Delta_\calM}}
\newcommand{\DelMo}{{\Delta_\calM^o}}
\newcommand{\calbM}{\bar{\cal M}}
\newcommand{\calF}{{\cal F}}
\newcommand{\calG}{{\cal G}}
\newcommand{\calB}{{\cal B}}
\newcommand{\calb}{{\cal b}}
\newcommand{\calD}{{\cal D}}
\newcommand{\calH}{{\cal H}}
\newcommand{\calL}{{\cal L}}
\newcommand{\calK}{{\cal K}}
\newcommand{\calC}{{\cal C}}
\newcommand{\calT}{{\cal T}}
\newcommand{\calZ}{{\cal Z}}
\newcommand{\calU}{{\cal U}}
\newcommand{\calE}{{\cal E}}
\newcommand{\calEc}{{\calE_c}}
\newcommand{\calN}{{\cal N}}
\newcommand{\calV}{{\cal V}}
\newcommand{\calR}{{\cal R}}
\newcommand{\calY}{{\cal Y}}
\newcommand{\calX}{{\cal X}}
\newcommand{\dcalV}{\dot{\cal V}}
\newcommand{\dcalU}{\dot{\cal U}}
\newcommand{\dcalE}{{\dot{\cal E}}}
\newcommand{\dbq}{{\dot{\bf q}}}
\newcommand{\bu}{{\bf u}}
\newcommand{\bara}{\bar{ a}}
\newcommand{\barba}{\bar{\bf a}}
\newcommand{\barbE}{\bar{\bf E}}
\newcommand{\barbC}{\bar{\bf C}}
\newcommand{\barbF}{\bar{\bf F}}
\newcommand{\barbM}{\bar{\bf M}}
\newcommand{\barbN}{\bar{\bf N}}
\newcommand{\barbS}{\bar{\bf S}}
\newcommand{\barbT}{\bar{\bf T}}
\newcommand{\barbn}{\bar{\bf n}}
\newcommand{\barbp}{\bar{\bf p}}
\newcommand{\barbq}{\bar{\bf q}}
\newcommand{\barbu}{\bar{\bf u}}
\newcommand{\barf}{\bar{f}}
\newcommand{\barbf}{\bar{\bf f}}
\newcommand{\barbr}{\bar{\bf r}}
\newcommand{\barbx}{\bar{\bf x}}
\newcommand{\barby}{\bar{\bf y}}
\newcommand{\barm}{{\bar{m}}}
\newcommand{\barw}{{\bar{w}}}
\newcommand{\barb}{{\bar{b}}}
\newcommand{\bare}{{\bar{e}}}
\newcommand{\barimath}{{\bar{\imath}}}
\newcommand{\barbb}{{\bar{\bf b}}}
\newcommand{\barS}{{\bar{S}}}
\newcommand{\bart}{{\bar{t}}}
\newcommand{\barbt}{{\bar{\bf t}}}
\newcommand{\barbs}{{\bar{\bf s}}}
\newcommand{\barbv}{{\bar{\bf v}}}
\newcommand{\barnabla}{\bar{\nabla}}
\newcommand{\bgra}{{\bf \nabla}}
\newcommand{\dW}{{\dot{W}}}
\newcommand{\dU}{\dot{U}}
\newcommand{\dlam}{{\dot{\lambda}}}
\newcommand{\barlam}{{\bar{\lambda}}}
\newcommand{\bs}{{\bf s}}
\newcommand{\bb}{{{\bf b}}}
\newcommand{\bd}{{\bf d}}
\newcommand{\bl}{{\bf l}}
\newcommand{\barl}{{\bar{\bf l}}}
\newcommand{\dotbS}{\dot{\bf S}}
\newcommand{\dotbu}{\dot{\bf u}}
\newcommand{\dotbp}{\dot{\bf p}}
\newcommand{\dotbx}{\dot{\bf x}}
\newcommand{\dE}{\dot{\bf E}}
\newcommand{\bP}{{\bf P}}
\newcommand{\bS}{{\bf S}}
\newcommand{\dtau}{\,\mbox{d}\tau}
\newcommand{\dI}{{\,\mbox{d}I}}
\newcommand{\dt}{{\,\mbox{d}t}}
\newcommand{\ds}{{\,\mbox{d}s}}
\newcommand{\dd}{\mbox{d}}
\newcommand{\ddt}{\frac{\d}{\dt}}
\newcommand{\dbt}{{\dot{\bf t}}}
\newcommand{\dv}{\dot{\bf v}}
\newcommand{\dotv}{\dot{v}}
\newcommand{\du}{\dot{\bf u}}
\newcommand{\dotu}{\dot{ u}}
\newcommand{\dotxi}{\dot{\xi}}
\newcommand{\be}{{\bf e}}
\newcommand{\de}{\dot{\bf e}}
\newcommand{\dS}{\mbox{ d}S}
\newcommand{\dx}{\mbox{ d}x}
\newcommand{\ddx}{\frac{\d}{\dx}}
\newcommand{\dy}{\,\mbox{d}y}
\newcommand{\dO}{\,\mbox{d}\Oo}
\newcommand{\dOt}{\,\mbox{d}\Oot}
\newcommand{\dG}{\,\mbox{d} \Gamma}
\newcommand{\dbM}{{\rm d} \partial\calM}
\newcommand{\dM}{{\rm d} \calM}
\newcommand{\bU}{{\bf U}}
\newcommand{ \grad }{{\mbox{grad}}}
\newcommand{\alp}{{\alpha}}
\newcommand{\ab}{{{\alpha\beta}}}
\newcommand{ \eps}{{\epsilon}}
\newcommand{ \Ups}{{\Upsilon}}
\newcommand{ \sig}{{\sigma}}
\newcommand{ \Lam}{{\Lambda}}
\newcommand{ \barLam}{\bar{\Lambda}}
\newcommand{ \Gam}{{\Gamma}}
\newcommand{ \lam}{{\lambda}}
\newcommand{\bgamma}{{\mbox{\boldmath$\gamma$}}}
\newcommand{ \xx}{{\bf x}}
\newcommand{ \Xx}{{\bf X}}
\newcommand{ \Ff }{{\bf F }}
\newcommand{\barR }{{\bar{\real}}}
\newcommand{\barOo }{{\bar{\Omega}}}
\newcommand{\bc}{{\bf c}}
\newcommand{\bD}{{\bf D}}
\newcommand{\bZ}{{\bf Z}}
\newcommand{\epi}{{\mbox{epi }}}
\newcommand{\dom}{{\mbox{dom }}}
\newcommand{\Int}{{\mbox{int}}}
\newcommand{\tr}{{\mbox{tr}}}
\newcommand{\Lin}{{\cal{M}}}
\newcommand{\Cof}{{\mbox{Cof }}}
\newcommand{\rank}{{\mbox{rank }}}
\newcommand{\Diag}{{\mbox{Diag }}}
\newcommand{\Ker}{{\mbox{Ker }}}
\newcommand{\lin}{{\mbox{lin}}}

\newcommand{\Pis }{\Pi_{\sup}}
\newcommand{\Pist }{\Pi_{\small{\sta}}}
\newcommand{\Pii }{\Pi_{\inf}}

\newcommand{\com}{{\mbox{com }}}
\newcommand{\ine}{{\mbox{ine }}}
\newcommand{\Vol}{{\mbox{Vol}}}
\newcommand{\Inv}{{\mbox{Inv}}}
\newcommand{\dive}{{\mbox{div}}}
\newcommand{\curl}{{\mbox{curl}}}
\newcommand{\Sym}{{\mbox{Sym}}}
\newcommand{\Global}{{\mbox{Global}}}
\newcommand{\sym}{{\mbox{sym}}}
\newcommand{\divM}{{\mbox{div}}_\calM}
\newcommand{\VtoR}{\calU \rightarrow \real}
\newcommand{\mboxin}{${\mbox{ in }}$}
\newcommand{\mboxon}{{\mbox{ on }}}
\newtheorem{assumption}{Assumption}
\newtheorem{remark}{Remark}
\newtheorem{algorithm}{Algorithm}
\newtheorem{thm}{Theorem}
\newtheorem{Corrolary}{Corrolary}
\newtheorem{definition}{Definition}
\newtheorem{problem}{Problem}
\newtheorem{example}{Example}

\newcommand\realmp{{\real^m_+}}
\newcommand\realn{{\real_-}}
\newcommand\realmn{{\real^m_-}}
\newcommand\realp{{\real_+}}
\renewcommand\pl{{|}}
\renewcommand\calPs{{\calP_{\sta}}}
\renewcommand\AA{{A}}
\renewcommand\bA{{{A}}}

\renewcommand\Pi{{{P}}}
\newcommand\BB{{B}}

\renewcommand\II{{\calI}}
\newcommand\FF{{F}}
\newcommand\yy{{y}}
\newcommand\byd{{y^*}}
\newcommand\KK{{K}}
\newcommand\y{{y}}

\newcommand\xxd{{{x^*}}}
\renewcommand\barbx{{\bar{x}}}
\newcommand\x{{{x}}}
\renewcommand\bff{{{f}}}
\renewcommand\bc{{{f}}}
\newcommand\yyd{{y^*}}
\newcommand\baryd{{\bar{y}^*}}
\newcommand\barWd{{\barW^*}}
\newcommand\Wv{{\check{W}}}
\newcommand\calYd{\calY^*}
\newcommand\calXd{\calX^*}
\newcommand\G{{G\^{a}teaux} }

\renewcommand\eb{\begin{equation}}
\renewcommand\ee{\end{equation}}

\newcommand\calPd{{\calP^d}}
\newcommand\bepsd{{\beps^*}}
\newcommand\rhod{{\rho^*}}
\newcommand\PPd{{P^d}}
\newcommand\barbepsd{{\bar{\beps}^*}}
\newcommand\barxm{{{\bar{\bx}_\mu}}}
\newcommand\barhod{{\bar{\rho}^*}}
\renewcommand\calXk{{{\calX_f}}}
\newcommand\barbnu{{\bar{\bnu}}}
\renewcommand\AA{{A}}
\newcommand\QQ{{Q}}
\newcommand\XI{{\Xi}}
\newcommand\DD{{D}}
\newcommand\RR{{R}}
\renewcommand\SS{{S}}
\newcommand\bq{\bf{q}}
\newcommand\barblam{\bar{\blam}}
\newcommand\QQd{{Q_{d}}}
\newcommand\WW{{W}}
\newcommand\VV{V}
\renewcommand\bA{{{A}}}
\renewcommand\bI{{{I}}}
\renewcommand\Pi{{{P}}}
\renewcommand\bb{{\bf b}}
\newcommand\GG{{G}}
\newcommand\GGs{{G}}
\newcommand\LLo{L_o}
\renewcommand\II{{\calI}}
\newcommand\xxo{{{\bf x}^o}}
\renewcommand\xx{{{\bf x}}}
\renewcommand\barbx{{\bar{\bf x}}}
\renewcommand\yy{{{\bf y}}}
\renewcommand\zz{{{\bf z}}}
\renewcommand\ww{{{\bf w}}}
\renewcommand\barww{{{\bar{\ww}}}}
\renewcommand\bc{{{\bf c}}}
\renewcommand\barbx{{\bar{\xx}}}
\newcommand\Nb{N_{A}}
\renewcommand\bff{{{\bf c}}}

\begin{center}
{\large \textbf{Canonical Solutions to
Nonconvex \\
Minimization Problems over Lorentz Cone}
\vspace{0.4cm}\\[0pt]}
{\textbf{Ning Ruan$^{1,2}$ and  David Yang Gao$^{1}$}}
\vspace{.3cm} \\
{\small {\it 1. Graduate School of Information Technology and Mathematical Science, \\
University of Ballarat, Ballarat, Vic 3353, Australia. \\
2. Department of Mathematics and  Statistics, \\
Curtin University of Technology, Perth,  WA 6845, Australia. }}
\end{center}

\begin{abstract}
This paper presents a canonical dual approach for solving
nonconvex quadratic minimization problem. By using the canonical
duality theory, nonconvex primal minimization problems over
n-dimensional Lorentz cone can be transformed into certain canonical
dual problems with only one dual variable, which can be solved by
using standard convex minimization methods. Extremality conditions
of these solutions are  classified by the triality theory.
Applications are illustrated.

\end{abstract}

\section{Primal Problem and It's Canonical Dual}
The primal problem $(\calP)$ proposed to solve is the so-called
second order cone programming:
\eb
(\calP): \;\;\;\; \min \left\{
P(\xx) =  \half \la \xx ,  \QQ\xx \ra  -   \la  \xx , \bc \ra \; :
\; \; \xx \in \calC \; \right\} , \label{eq-ncp}
\ee
where $\QQ \in
\real^{n\times n} $ is a given symmetrical matrix; $\bc \in \real^n$
is a given vector; $\calC \subset \real^n$ is the so-called Lorentz cone in $\real^n$:
\eb
\calC = \{ (x_1, \xx_2) \in \real \times \real^{n-1} | \;\; \|  \xx_2 \| \le  x_1,~ x_1> 0\},
\ee
which is a special second order cone.
The problem $(\calP)$ appears in many applications such as
structural optimization, filter design, and grasping force
optimization in robotics. Extensive research has been focused on this subject.

In this paper we present a canonical dual approach for solving the
second order cone optimization problem $(\calP)$. By using the
canonical duality theory developed in
\cite{gao-book00,gao-jogo00,gao-opt03, gao-jogo04}, the canonical
dual problem $(\calPd)$  can be formulated as
\eb
\max \left\{ \PPd(\sig) = - \half \la \bc , \; \GG^{-1}(\sig) \bc \ra\; : \;\; \sig  \in \calS_a \right\} ,
\ee
where
$\GG (\sig) = \QQ + \sig \LLo $,
\[
\LLo = \left( \begin{array}{cc}
-1 & 0 \\
0 & \II_{n-1} \end{array}
\right)
\]
is the Lorentz matrix, in which, $\II_{n-1}$ is an identical matrix in $\real^{(n-1)\times (n-1)}$.
The dual feasible space $\calS_a \subset \real$ is defined by
\eb
\calS_a = \{ \sig \in \real\; | \;
\sig\ge  0, \;\; \det \GG  (\sig) \neq 0 \} .
\ee
\begin{thm}\label{primal solu}
The problem $(\calPd)$ is canonically dual to $(\calP)$ in the sense that
the vector $\barsig \in \calS_a$ is a KKT point of $(\calPd)$  if and only if  the vector
\eb
\barbx = \GG^{-1} (\barsig) \bc
\ee
is a KKT point  of  $(\calP)$, and
\eb
\PP(\barbx) = \PP^d(\barsig).
\ee
\end{thm}

\noindent
{\em Proof}.
By the standard procedure of the canonical dual transformation, we rewrite the  cone  constraint
$\| \bx_2\| \le x_1$ in the quadratic form $\half \bx^T \LLo \bx \le 0  $ and
introduce  a nonlinear transformation (i.e. the so-called {\em geometrical mapping})
$\eps = \Lam(\xx) = \half \bx^T \LLo \bx: \real^n \rightarrow \real $.
Thus, the cone constraint  $\bx \in \calC$ can be replaced identically by $\eps(\xx) = \Lam(\xx) \le 0$. Let
\eb
\VV(\eps) = \left\{ \begin{array}{ll}
0 & \mbox{ if } \eps  \le  {0}  , \\
+ \infty & \mbox{ otherwise}.
\end{array}
\right. \label{eq-wws}
\ee
The  primal problem $(\calP)$ can be written in the following canonical form \cite{gao-opt03}:
\eb
(\calP_c): \;\; \min \left\{ \PP(\xx) = \VV(\Lam(\xx)) + \half \xx^T \QQ \xx -\bff^T \xx   \; | \;\;
\xx \in \real^n  \right\} . \label{eq-cpp}
\ee
According to the Fenchel  transformation, the sup-conjugate $\VV^\sharp$
of the function $\VV(\eps)$ is defined by
\[
\VV^\sharp(\sig) = \sup_{\eps \in \real} \{ \eps^T \sig - \VV(\eps) \} =
\left\{ \begin{array}{ll} 0 & \mbox{ if } \sig \ge  {0}, \\
+ \infty & \mbox{ otherwise}.
\end{array}
\right.
\]
Since $\VV(\eps)$ is a proper closed convex function over $\real_- := \{ \eps
\in \real | \; \eps \le {  0}\}$, we know that
\eb
\sig \in \partial \VV(\eps)  \;\; \Leftrightarrow  \;\;\; \eps \in \partial \VV^\sharp(\sig)
\;\; \Leftrightarrow  \;\;\; \VV(\eps) +  \VV^\sharp(\sig) = \eps  \sig .\label{eq-canodual}
\ee
The pair of $(\eps, \sig)$ is then
called a {\em generalized canonical dual pair} on $\real_- \times \real_+$
by the definition introduced in \cite{gao-book00,gao-jogo00}.
Following the original idea of Gao and Strang \cite{gao-strang89}, we replace
$\VV(\Lam(\xx))$ in equation (\ref{eq-cpp}) by the Fenchel-Young equality
$\VV(\Lam(\xx)) = \Lam(\xx)^T \sig - \VV^\sharp(\sig)$. Then
the  so-called  {\em total complementary function}
$\Xi(\xx, \sig):\real^n \times \real_+ \rightarrow \real$
associated with the problem $(\calP_c)$ can be defined as below
\begin{eqnarray}
\Xi(\xx, \sig) = \Lam(\xx)^T \sig - \VV^\sharp (\sig)+ \half \xx^T \QQ \xx -\bff^T \xx .\label{eq-tcomf}
\end{eqnarray}
By the definition of $\Lam(\xx)$ and $\VV^\sharp (\sig)$, on $\real^n \times \real_+$ we have
\begin{eqnarray}
\Xi(\xx, \sig) =\half \xx^T \GGs(\sig) \xx - \xx^T  \bff  . \label{eq-llg}
\end{eqnarray}
The criticality condition of $\Xi(\xx, \sig)$ leads to the equilibrium equation
\eb
\GGs(\sig)   \xx =  \bff,\label{equili}
\ee
and the KKT conditions
\eb
\sig \ge 0 , \;\;  \xx^T \LLo \xx  \le 0 , \;\; \sig   \xx^T \LLo \xx    = 0.
\ee

Substituting (\ref{equili}) into the canonical dual transformation
\[
\PP^d(\sig) = \sta \{ \Xi(\xx, \sig) \; | \;\; \xx \in \calC \},
\]
the canonical dual function $\PP^d(\sig)$ is then formulated.

It is easy to prove that if  $\barsig \ge  0 $ is a KKT point of $(\calP^d )$, then we have
\begin{eqnarray}
&& \barsig \ge 0 , \;\; \;\;  \nabla\PPd(\barblam) =
\frac{1}{2}\barbx(\barbsig)^T L_0  \barbx(\barbsig)\le 0, \label{comprho} \\
&& \barsig\cdot(\frac{1}{2}\barbx(\barbsig)^T L_0  \barbx(\barbsig))=0, \label{compeps}
\end{eqnarray}
where $\barbx (\barsig) = G^{-1}(\barsig) \bff $.
This shows that $\barbx(\barsig)$ is also a KKT point of the primal problem $(\calP )$.

By the complementarity condition (\ref{compeps}),
and the fact of  $\barbx = G^{-1}(\sig) \bff $, we have
\begin{eqnarray*}
P^d(\barsig)
&=& -\frac{1}{2}\bc^T G^{-1}(\barsig) \bc\nonumber\\
&=& \frac{1}{2}\bc^T G^{-1}(\barsig) \bc-\bc^T G^{-1}(\barsig) \bc\nonumber\\
&=& \frac{1}{2}(G^{-1}(\barsig)\bc)^T G(\barsig)G^{-1}(\barsig) \bc
-\bc^T G^{-1}(\barsig) \bc\nonumber\\
&=& \frac{1}{2}(G^{-1}(\barsig)\bc)^T (Q+\barsig L_0)G^{-1}(\barsig)\bc
-\bc^T G^{-1}(\barsig) \bc\nonumber\\
&=& \frac{1}{2}(G^{-1}(\barsig)\bc)^T Q G^{-1}(\barsig)\bc-\bc^T G^{-1}(\barsig) \bc+
\frac{\barsig}{2}(G^{-1}(\barsig)\bc)^T L_0 G^{-1}(\barsig)\bc\nonumber\\
&=& \frac{1}{2}(G^{-1}(\barsig)\bc)^T Q G^{-1}(\barsig)\bc- \bc^T G^{-1}(\barsig) \bc\nonumber\\
&=& \frac{1}{2}\barbx^T Q \barbx -\bc^T\barbx\nonumber\\
&=& P(\barbx)
\end{eqnarray*}
This proves the theorem.\hfill $\Box$

\section{Extremality Conditions}
Let
\begin{eqnarray}
\calS_a^+ &=& \{ \sig \in \real\; | \;\sig \ge  0, \;\;\GG(\sig) \succ 0\}.
\end{eqnarray}

\begin{thm}\label{condition}
Suppose that $\barsig \in \calS_a$ is a solution to $(\calP)$ and
\[
\barbx = \GG^{-1} (\barsig)\bc .
\]

If $\barsig \in \calS_a^+$, then $\barbx$ is
a global minimizer of $\PP(\bx)$ on $\calC$ and
\eb
\PP(\barbx) =
\min_{\bx \in \calC} \PP(\bx) =  \max_{\sig  \in \calS_a^+}
\PP^d(\sig) = \PP^d(\barsig). \label{def}
\ee

\end{thm}

{\bf Proof.}
By Theorem \ref{primal solu} we know that  the  vector  $\barsig  \in \calS_a$
is a KKT point of the problem $(\calP^d)$ if and only if
$\barbx = \GG^{-1} (\barsig)\bc$ is a KKT point of the problem $(\calP)$, and
\eb
\PP(\barbx ) = \Xi(\barbx , \barsig) =
\PP^d(\barsig). \label{eq:triality}
\ee
Particularly, if $\barsig \in \calS_a^+$,  the canonical dual function  $\PP^d(\sig)$ is concave.
In this case, the total complementary function  $\Xi$  is a saddle function, i.e.,
it is convex in $\xx \in \real^n$ and concave in $\sig\in \calS^+_a$. Thus, we have
\begin{eqnarray*}
\PP^d(\barsig)& =& \max_{ \sig \in \calS_a^+}  \PP^d(\sig)\\
&=&  \max_{\sig\in\calS^+_a} \min_{\xx \in
\real^n} \Xi(\xx, \sig) = \min_{\xx \in \real^n} \max_{\sig \in\calS_a^+ }  \Xi(\xx, \sig) \\
&=&\min_{\xx \in \real^n}\left\{\half \xx^T   \QQ \xx -\xx^T \bff  + \max_{\sig \in\calS_a^+}
\left\{\left(\frac{1}{2}\xx^T L_0 \xx \right) \sig - \VV^\sharp (\sig)\right\} \right\}
= \min_{\xx \in \calC}  \PP(\bx)
\end{eqnarray*}
due to the fact that
\[
\VV(\Lam( \bx))= \max_{ \sig\in\calS_a^+ }
\left \{ \left(\frac{1}{2}\xx^T L_0 \xx \right) \sig- \VV^\sharp (\sig)\right \} =
\left\{ \begin{array}{ll}0 & \mbox{ if } \bx \in \calC, \\
\infty & \mbox{ otherwise}.
\end{array} \right.
\]
From Theorem \ref{primal solu}  we have  (\ref{def}).
\hfill $\Box$

In a special case when  $\QQ = \Diag (\bq)$ is a diagonal matrix
with $\bq = \{ q_i \} \in \real^n$ being its  diagonal elements, we  have \eb
G^{-1}(\sigma)= \left\{ \frac{1}{q_i+\sigma \delta_i^{-}} \right\},
\ee In this case,
\eb
P^d(\sigma)=-\frac{1}{2} \sum_{i=1}^n
\frac{c_i^2}{q_i+\sigma \delta_i^{-}},
\ee
where
\eb
\delta_i^{-}=\bigg\{\begin{array}{l l}
-1 &i=1\\
1&i\neq 1
\end{array},\;\;
\delta_i^{+}=\bigg\{\begin{array}{l l}
1 &i=1\\
0&i\neq 1\end{array}.
\ee

For the given $\{c_i\} $, and  $\{ q_i\}$ such that $-q_n\le q_{n-1}
\le \dots, \le q_1$, the dual variable $\sig$ can be solved
completely within each interval $-q_{i+1} < \sig < - q_i $ or
$-q_2< \sigma<q_1$, such that $q_i < q_{i+1}$ ($i=2, \cdots, n$).

\section{Applications}
We now list a few examples to illustrate the applications of the
theory presented in this paper.

\subsection{Two-D nonconvex minimization }
First of all, let us consider two dimensional concave minimization problem:
\eb
(\calP): \;\; \min \left\{ P(\xx) =  \frac{1}{2}(q_1
x_1^2+q_2 x_2^2)-c_1 x_1-c_2 x_2\; : \; \|  \xx_2 \| \le x_1,(x_1,
\xx_2) \in \real^2 \; \right\} ,
\ee

On the dual feasible set
\eb
\calS_a = \{ \sig \in \real\; | \;
\sig\ge  0, \;\;(q_{1}-\sigma)(q_{2}+\sigma) \neq 0\},
\ee
the canonical dual function has the form of
\eb
P^d(\sig )=-\frac{1}{2}[c_{1},c_{2}]^T
\left[\matrix{\frac{1}{q_{1}-\sig}&\cr &\frac{1}{q_{2}+\sig}}
\right] \left[\matrix{c_{1}\cr c_{2}} \right].
\ee

\begin{figure}[h!]
\centering
\mbox{\resizebox{!}{2.0in}{\includegraphics{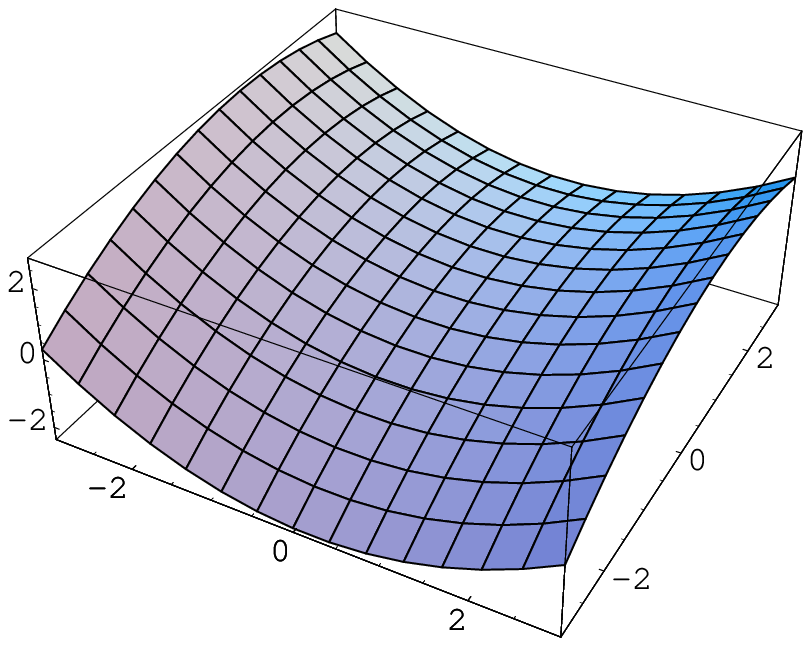}}\quad\quad
\resizebox{!}{2.0in}{\includegraphics{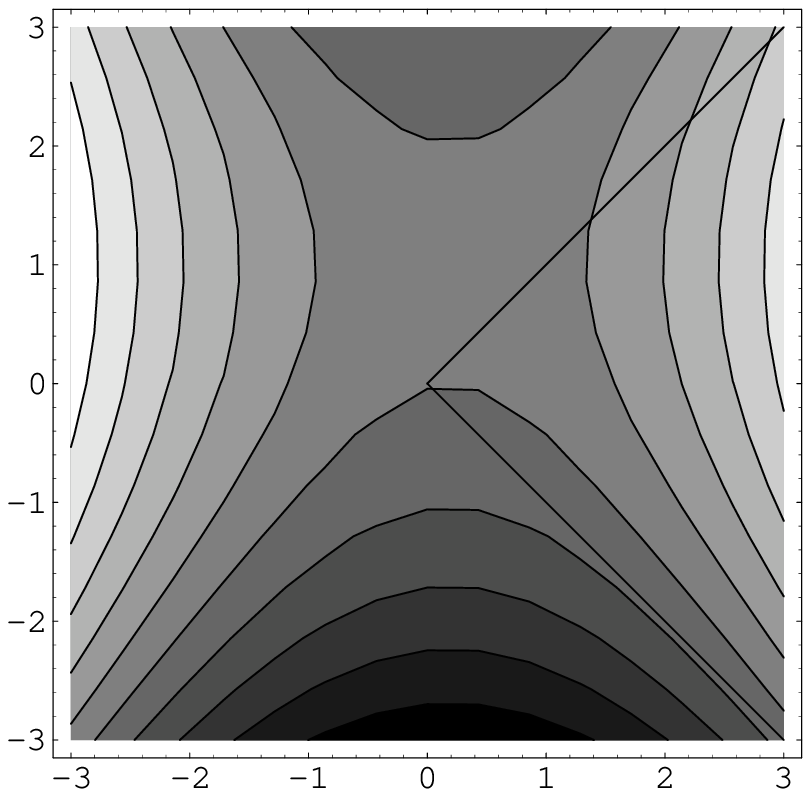}}}\vspace{-0.5cm}
\caption{\label{primal1} Graphs of $P(x)$ and its contour for two
dimensional problem.}
\end{figure}

\begin{figure}[h!]
\vspace{0.1cm} {\centering
\resizebox{!}{2.0in}{\includegraphics{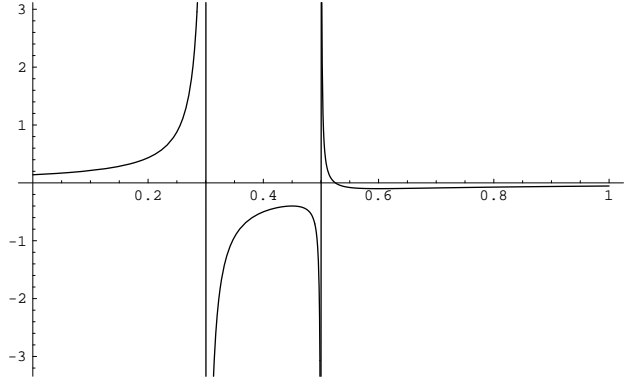}} \vspace{-0.5cm}
\caption{\label{dual1} Graph of $P^d(\sig)$ for two dimensional
problem.}}
\end{figure}

Assume $ q_1=0.1$, $q_2=-0.3$, $c_1=0.5$, $c_2=-0.3$, so we have
\eb
\sig=0.45\in \calS_a^+=\{\sig\in {\mathbb
R}|\;\; 0.3<\sig<0.5\}.
\ee
By Theorem \ref{primal solu}, we know that $x$ = $\{c_1/(q_1-\sig),
c_2/(q_2+\sig)\}$= $\{2,-2\}$ is a global minimizer, it is easy to
verify that $ P(x)=P^d(\sig)=-0.4 $(see Figure \ref{primal1}-\ref{dual1}).

\subsection{Two-D general nonconvex minimization }
\eb
(\calP): \;\; \min \left\{ P(\xx) =  \frac{1}{2}(q_1 x_1^2+q_2
x_2^2+2 q_3 x_1 x_2)-c_1 x_1-c_2 x_2\; : \; \|  \xx_2 \| \le
x_1,(x_1, \xx_2) \in \real^2 \; \right\} ,
\ee

On the dual feasible set
\eb
\calS_a = \{ \sig \in \real^2\; | \;
\sig\ge  0, \;\;(q_{1}-\sig)(q_{2}+\sig) -q_3^2\neq 0\},
\ee
the canonical dual function has the form of
\eb
P^d(\sig)=-\frac{1}{2}[c_{1},c_{2}]^T
\left[\matrix{q_1-\sig&q_3\cr a_3&q_2+\sig }\right]^{-1}
\left[\matrix{c_{1}\cr c_{2}} \right]
\ee

\begin{figure}[h!]
\centering
\mbox{\resizebox{!}{2.0in}{\includegraphics{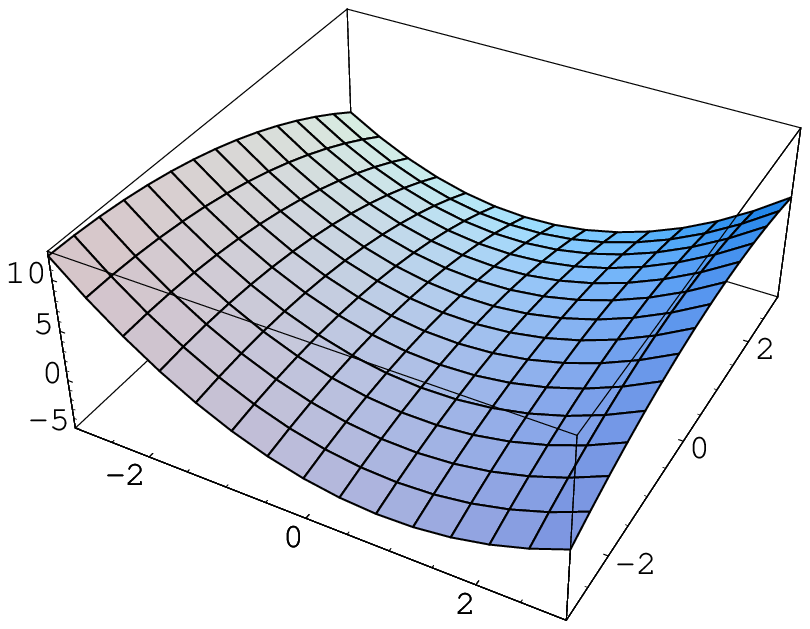}}\quad\quad
\resizebox{!}{2.0in}{\includegraphics{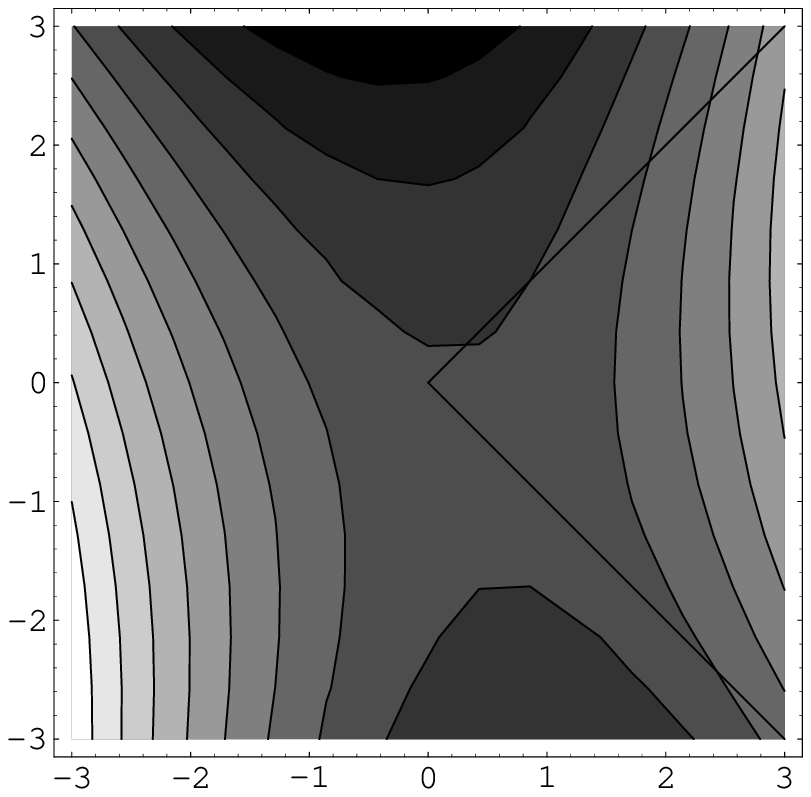}}}\vspace{-0.5cm}
\caption{\label{primal2} Graphs of $P(x)$ and its contour for
general two dimensional problem.}
\end{figure}

\begin{figure}[h!]
\vspace{0.1cm} {\centering
\resizebox{!}{2.0in}{\includegraphics{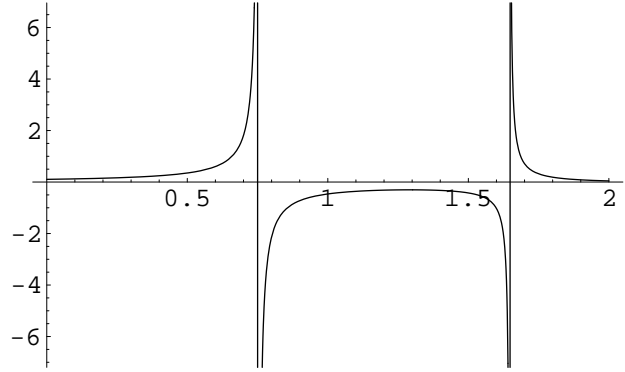}} \vspace{-0.5cm}
\caption{\label{dual2} Graph of $P^d(\sig)$ for general two
dimensional problem.}}
\end{figure}

If we choose $ q_1=1.8$, $q_2=-0.6$, $q_3=0.4$, $c_1=0.5$,
$c_2=0.6$, then we have
\eb
\sig=1.29\in \calS_a^+=\{\sigma\in {\mathbb R}|\;\;
(q_{1}-\sig)(q_{2}-\sig) -a_3^2>0\}.
\ee
By Theorem \ref{primal solu}, we know that $x$ = $\left[\matrix{q_1-\sig
&q_3\cr q_3&q_2+\sig}\right]^{-1} \left[\matrix{c_{1}\cr
c_{2}} \right]$= $\{0.5500,0.5499\}$ is a global minimizer, and
$P(x)=P^d(\sig)=-0.3025$(see Figure \ref{primal2}-\ref{dual2}).

\subsection{Three-D general nonconvex minimization }
\eb
(\calP): \;\; \min \left\{ P(\xx) =  \half \la \xx ,  \QQ\xx \ra
- \la  \xx , \cc \ra\; : \; \| \xx_2 \| \le x_1,(x_1, \xx_2) \in \real^3 \; \right\} ,
\ee
where
$Q=\left[\matrix{q_{11}&q_{12}&q_{13}\cr q_{12}&q_{22}&q_{23}\cr
q_{13}&q_{23}&q_{33}}\right]$, $c=\{c_1,c_2,c_3\}$.  On the dual feasible set
\eb
\calS_a = \{ \sig \in \real\; | \;
\sig  \ge 0 , \;\;{\rm det} G(\sig) \neq 0 \},
\ee
where
\eb
G(\sig) =\left[\matrix{q_{11}-\sig &q_{12}&q_{13}\cr
q_{12}&q_{22}+\sig&q_{23}\cr q_{13}&q_{23}&q_{33}+\sig }\right],
\ee
the canonical dual function has the form of
\eb
P^d(\sigma)=-\frac{1}{2}[c_{1},c_{2},c_3]^T G^{-1}(\sig)
\left[\matrix{c_{1}\cr c_{2}\cr c_3} \right]
\ee

\begin{figure}[h!]
\vspace{0.1cm} {\centering
\resizebox{!}{2.0in}{\includegraphics{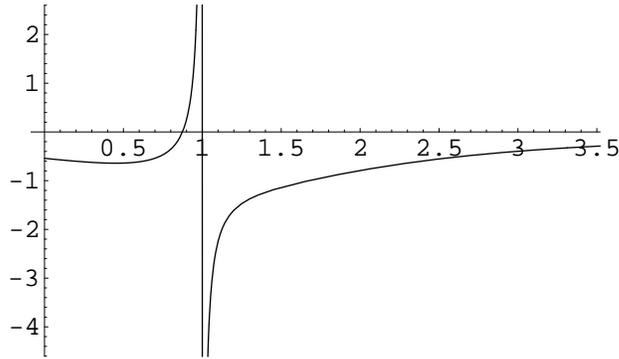}} \vspace{-0.5cm}
\caption{\label{dual3} Graph of $P^d(\sig)$ for general three
dimensional problem.}}
\end{figure}

Suppose $ q_{11}=2$, $q_{22}=-2$, $q_{33}=1$, $q_{12}=-1$,
$q_{13}=2$, $q_{23}=0$, $c_1=1.5$, $c_2=-0.5$, $c_3=1.5$, we have
\eb
\sigma=0.4509\in \calS_a^+=\{\sigma\in {\mathbb R}|\;\; \sigma \geq 0, G(\sigma)\succ 0\}.
\ee
By theorem \ref{primal solu}, we know that $x$ = $ G^{-1}(\sig)
\left[\matrix{ c_{1}\cr c_{2}\cr c_3} \right]$=
$\{0.4355,0.0416,0.4335\}$ is a global minimizer, and $P(x)=P^d(\sig)=-0.6413$(see Figure \ref{dual3}).

\section{Conclusions}
We have presented a concrete application of the canonical dual transformation
and triality to conic optimization problems. Results show that by the use of this method,
the nonconvex cone constrained problem $(\calP)$ in $\real^n$ can be reformulated as a perfect
dual problem in $\real$, also the $KKT$ points and extremality conditions of the originally difficult
problems are idetified by Theorem \ref{primal solu} and \ref{condition}.  Physically speaking,
each optimal point represents a stable equilibrium state of the system. Duality theory reveals
the intrinsic pattern of duality relations of these critical points, and plays an important role
in nonconvex analysis, detailed study and comprehensive applications of this theory were presented in
monograph \cite{gao-book00}.

\section*{Acknowledgment}
 This paper was partially supported by a grant (AFOSR FA9550-10-1-0487)
from the US Air Force Office of Scientific Research. Dr. Ning Ruan was
supported by a funding from the Australian Government
under the Collaborative Research Networks (CRN) program.

\end{document}